\documentclass{article}
\usepackage{amssymb,amsmath,theorem,euscript}

\input{epsf}

\newcounter{sec}

\def\sm{\smallskip}

\newcounter{punct}[sec]

\def\punct{\refstepcounter{punct}{\arabic{sec}.\arabic{punct}.  }} 

\def\COUNTERS{\addtocounter{sec}{1}
              \setcounter{punct}{0}
          \setcounter{equation}{0}
          \setcounter{theorem}{0}
                  }

\newtheorem{theorem}{Theorem}[sec]
\newtheorem{proposition}[theorem]{Proposition}
\newtheorem{lemma}[theorem]{Lemma}
\newtheorem{corollary}[theorem]{Corollary}

\begin{document}

 \def\ov{\overline}
\def\wt{\widetilde}
\def\wh{\widehat}
 \newcommand{\rk}{\mathop {\mathrm {rk}}\nolimits}
\newcommand{\Aut}{\mathop {\mathrm {Aut}}\nolimits}
\newcommand{\Out}{\mathop {\mathrm {Out}}\nolimits}
 \newcommand{\tr}{\mathop {\mathrm {tr}}\nolimits}
  \newcommand{\diag}{\mathop {\mathrm {diag}}\nolimits}
  \newcommand{\supp}{\mathop {\mathrm {supp}}\nolimits}
  \newcommand{\indef}{\mathop {\mathrm {indef}}\nolimits}
  \newcommand{\dom}{\mathop {\mathrm {dom}}\nolimits}
  \newcommand{\im}{\mathop {\mathrm {im}}\nolimits}
   \newcommand{\ind}{\mathop {\mathrm {ind}}\nolimits}
    \newcommand{\codim}{\mathop {\mathrm {codim}}\nolimits}
 
\renewcommand{\Re}{\mathop {\mathrm {Re}}\nolimits}

\def\Br{\mathrm {Br}}

\def\SL{\mathrm {SL}}
\def\Diag{\mathrm {Diag}}
\def\SU{\mathrm {SU}}
\def\GL{\mathrm {GL}}
\def\U{\mathrm U}
\def\OO{\mathrm O}
 \def\Sp{\mathrm {Sp}}
 \def\SO{\mathrm {SO}}
\def\SOS{\mathrm {SO}^*}
 \def\Diff{\mathrm{Diff}}
 \def\Vect{\mathfrak{Vect}}
\def\PGL{\mathrm {PGL}}
\def\PU{\mathrm {PU}}
\def\PSL{\mathrm {PSL}}
\def\Symp{\mathrm{Symp}}
\def\End{\mathrm{End}}
\def\Mor{\mathrm{Mor}}
\def\Aut{\mathrm{Aut}}
 \def\PB{\mathrm{PB}}
 \def\cA{\mathcal A}
\def\cB{\mathcal B}
\def\cC{\mathcal C}
\def\cD{\mathcal D}
\def\cE{\mathcal E}
\def\cF{\mathcal F}
\def\cG{\mathcal G}
\def\cH{\mathcal H}
\def\cJ{\mathcal J}
\def\cI{\mathcal I}
\def\cK{\mathcal K}
 \def\cL{\mathcal L}
\def\cM{\mathcal M}
\def\cN{\mathcal N}
 \def\cO{\mathcal O}
\def\cP{\mathcal P}
\def\cQ{\mathcal Q}
\def\cR{\mathcal R}
\def\cS{\mathcal S}
\def\cT{\mathcal T}
\def\cU{\mathcal U}
\def\cV{\mathcal V}
 \def\cW{\mathcal W}
\def\cX{\mathcal X}
 \def\cY{\mathcal Y}
 \def\cZ{\mathcal Z}
\def\0{{\ov 0}}
 \def\1{{\ov 1}}
 \def\frA{\mathfrak A}
 \def\frB{\mathfrak B}
\def\frC{\mathfrak C}
\def\frD{\mathfrak D}
\def\frE{\mathfrak E}
\def\frF{\mathfrak F}
\def\frG{\mathfrak G}
\def\frH{\mathfrak H}
\def\frI{\mathfrak I}
 \def\frJ{\mathfrak J}
 \def\frK{\mathfrak K}
 \def\frL{\mathfrak L}
\def\frM{\mathfrak M}
 \def\frN{\mathfrak N} \def\frO{\mathfrak O} \def\frP{\mathfrak P} \def\frQ{\mathfrak Q} \def\frR{\mathfrak R}
 \def\frS{\mathfrak S} \def\frT{\mathfrak T} \def\frU{\mathfrak U} \def\frV{\mathfrak V} \def\frW{\mathfrak W}
 \def\frX{\mathfrak X} \def\frY{\mathfrak Y} \def\frZ{\mathfrak Z} \def\fra{\mathfrak a} \def\frb{\mathfrak b}
 \def\frc{\mathfrak c} \def\frd{\mathfrak d} \def\fre{\mathfrak e} \def\frf{\mathfrak f} \def\frg{\mathfrak g}
 \def\frh{\mathfrak h} \def\fri{\mathfrak i} \def\frj{\mathfrak j} \def\frk{\mathfrak k} \def\frl{\mathfrak l}
 \def\frm{\mathfrak m} \def\frn{\mathfrak n} \def\fro{\mathfrak o} \def\frp{\mathfrak p} \def\frq{\mathfrak q}
 \def\frr{\mathfrak r} \def\frs{\mathfrak s} \def\frt{\mathfrak t} \def\fru{\mathfrak u} \def\frv{\mathfrak v}
 \def\frw{\mathfrak w} \def\frx{\mathfrak x} \def\fry{\mathfrak y} \def\frz{\mathfrak z} \def\frsp{\mathfrak{sp}}
 \def\bfa{\mathbf a} \def\bfb{\mathbf b} \def\bfc{\mathbf c} \def\bfd{\mathbf d} \def\bfe{\mathbf e} \def\bff{\mathbf f}
 \def\bfg{\mathbf g} \def\bfh{\mathbf h} \def\bfi{\mathbf i} \def\bfj{\mathbf j} \def\bfk{\mathbf k} \def\bfl{\mathbf l}
 \def\bfm{\mathbf m} \def\bfn{\mathbf n} \def\bfo{\mathbf o} \def\bfp{\mathbf p} \def\bfq{\mathbf q} \def\bfr{\mathbf r}
 \def\bfs{\mathbf s} \def\bft{\mathbf t} \def\bfu{\mathbf u} \def\bfv{\mathbf v} \def\bfw{\mathbf w} \def\bfx{\mathbf x}
 \def\bfy{\mathbf y} \def\bfz{\mathbf z} \def\bfA{\mathbf A} \def\bfB{\mathbf B} \def\bfC{\mathbf C} \def\bfD{\mathbf D}
 \def\bfE{\mathbf E} \def\bfF{\mathbf F} \def\bfG{\mathbf G} \def\bfH{\mathbf H} \def\bfI{\mathbf I} \def\bfJ{\mathbf J}
 \def\bfK{\mathbf K} \def\bfL{\mathbf L} \def\bfM{\mathbf M} \def\bfN{\mathbf N} \def\bfO{\mathbf O} \def\bfP{\mathbf P}
 \def\bfQ{\mathbf Q} \def\bfR{\mathbf R} \def\bfS{\mathbf S} \def\bfT{\mathbf T} \def\bfU{\mathbf U} \def\bfV{\mathbf V}
 \def\bfW{\mathbf W} \def\bfX{\mathbf X} \def\bfY{\mathbf Y} \def\bfZ{\mathbf Z} \def\bfw{\mathbf w}
 \def\R {{\mathbb R }} \def\C {{\mathbb C }} \def\Z{{\mathbb Z}} \def\H{{\mathbb H}} \def\K{{\mathbb K}}
 \def\N{{\mathbb N}} \def\Q{{\mathbb Q}} \def\A{{\mathbb A}} \def\T{\mathbb T} \def\P{\mathbb P} \def\G{\mathbb G}
 \def\bbA{\mathbb A} \def\bbB{\mathbb B} \def\bbD{\mathbb D} \def\bbE{\mathbb E} \def\bbF{\mathbb F} \def\bbG{\mathbb G}
 \def\bbI{\mathbb I} \def\bbJ{\mathbb J} \def\bbK{\mathbb K} \def\bbL{\mathbb L} \def\bbM{\mathbb M} \def\bbN{\mathbb N} \def\bbO{\mathbb O}
 \def\bbP{\mathbb P} \def\bbQ{\mathbb Q} \def\bbS{\mathbb S} \def\bbT{\mathbb T} \def\bbU{\mathbb U} \def\bbV{\mathbb V}
 \def\bbW{\mathbb W} \def\bbX{\mathbb X} \def\bbY{\mathbb Y} \def\kappa{\varkappa} \def\epsilon{\varepsilon}
 \def\phi{\varphi} \def\le{\leqslant} \def\ge{\geqslant}

\def\UU{\bbU}
\def\Mat{\mathrm{Mat}}
\def\tto{\rightrightarrows}

\def\Gr{\mathrm{Gr}}

\def\graph{{graph}}

\def\O{\mathrm{O}}

\def\la{\langle}
\def\ra{\rangle}

\def\B{\mathrm B}
\def\Int{\mathrm{Int}}
\def\LGr{\mathrm{LGr}}


\def\I{\mathbb I}
\def\M{\mathbb M}
\def\T{\mathbb T}
\def\S{\mathrm S}

\def\Lat{\mathrm{Lat}}
\def\LLat{\mathrm{LLat}} 
\def\Mod{\mathrm{Mod}}
\def\LMod{\mathrm{LMod}}
\def\Naz{\mathrm{Naz}}
\def\naz{\mathrm{naz}}
\def\bNaz{\mathbf{Naz}}
\def\AMod{\mathrm{AMod}}
\def\ALat{\mathrm{ALat}}
\def\MAT{\mathrm{MAT}}
\def\Mar{\mathrm{Mar}}

\def\Ver{\mathrm{Vert}}
\def\Bd{\mathrm{Bd}}
\def\We{\mathrm{We}}
\def\Heis{\mathrm{Heis}}
\def\Pol{\mathrm{Pol}}
\def\Ams{\mathrm{Ams}}
\def\Herm{\mathrm{Herm}}

\def\F{\bbF}
\def\V{\bbV}

\def\Exp{\mathsf{Exp}}

 \newcommand{\gr}{\mathop {\mathrm {gr}}\nolimits}

\begin{center}
 \Large \bf
 
 On the Weil representation of infinite-dimensional symplectic group over a finite field

 \large \sc
 
 \bigskip
 
 Yu.A.Neretin%
 \footnote{Supported by the grant FWF, Project P28421.}
 
\end{center}

\bigskip

{\small We extend the Weil representation of infinite-dimensional symplectic group
	 to a representation a certain category
	of linear relations.}

\section{The statement of the paper}

\COUNTERS

{\bf \punct The spaces $\V_{2 \mu}$.} 
Denote by $\C^\times$ be the multiplicative
group of complex numbers.
 Denote by $\F$ a finite field 
of characteristic $p>2$, by $\F^\times$ its multiplicative group.
Fix a nontrivial character $\Exp(\cdot)$ of the additive group of $\F$,
$$
\Exp(a+b)=\Exp(a)\,\Exp(b), \qquad \Exp(\cdot)\in \C^\times .
$$
Values of $\Exp(\cdot)$ are $e^{2\pi i k/p}$.

For $n=0$, $1$, \dots,  we denote by $V_n^d$ an $n$-dimensional linear space over $\F$,
to be definite we assume that  $V_n^d$ is the coordinate space $\F^n$. Denote by $V_n^c$ the second copy of this space.

Next, we define the space $V_\infty^d$ as a direct sum of a countable number of $\F$
and the space $V_\infty^c$ as a direct product of a countable number of $\F$.
We equip the first space by the discrete topology, the second space by 
the product-topology ('$d$' is an abreviation of 'discrete' and '$c$' of 'compact). 

Below subscripts $\mu$, $\nu$, $\kappa$ denote $0$, $1$, \dots, $\infty$.

A space $V^c_\mu$ is the Pontryagin dual of
$V^d_\mu$, the pairing is given by
$$
[x,y]:=\sum_{j=1}^\mu x_j y_j.
$$
(actually, the sum is finite).
For a subspace $L\subset V^d_\mu$ we denote by $L^\circ$ the {\it anihilator}
of $L$ in $V^c_\mu$ and vice versa.

Denote by $\V_{2\mu}$ the space
$$
\V_{2\mu}:=V_\mu^d\oplus V_\mu^c
$$
equipped with the following skew-symmetric bilinear form ({\it symplectic form}):
$$
\bigl\{(x,y), (x',y')\bigr\}=\sum_{j=1}^\mu (x_j y_j'- y_j x_j').
$$
The space $\V_{2\infty}$ is a locally compact Abelian group and we can apply
the Pontryagin duality.

For a subspace $R\subset \V_{2\mu}$ we denote by $R^\lozenge$ the {\it orthocomplement} with respect
to $\{\cdot,\cdot\}$. By the Pontryagin duality (see \cite{Mor}, Proposition 38),
for closed subspaces $X$, $Y\subset V_{2\infty}$ the conditions
$X=Y^\lozenge$ and $Y=X^\lozenge$ are equivalent (in other words for a closed
subspace $X$ we have $X^{\lozenge\lozenge}=X$).

Denote by $\Sp(\V_{2\mu})$ the group of invertible continuous operators 
in $\V_{2\mu}$
preserving the form $\{\cdot,\cdot\}$. 

\sm

{\bf\punct The Weil representation of $\Sp(\V_{2\infty})$.}
Consider the space $\ell^2(V^d_{\mu})$ on the set  $V_{\mu}^d$.
For any $v=(v^d,v^c)\in \V_{2\mu}$ we define a unitary operator $a(v)$ in $\ell^2(V^d_{\mu})$
given by
$$
a(v)\,f(x)=f(x+v^d)\,\Exp\Bigl[
 \sum \bigl( x_j v_j^c+ \tfrac 12 v_j^d v_j^c\bigr)  \Bigr].
$$
We have
$$
a(v) a(w)= \Exp\bigl( \tfrac 12 \{v,w \} \bigr) \, a(v+w).
$$
Thus we get a  projective unitary representation of the additive group $\V_{2\mu}$ 
or a linear representation of the {\it Heisenberg group}, the latter group 
is $\V_{2\mu}\oplus \F$ with a
product
$$
(v;s)\star  (w;t)=\bigl(v+w; s+t+\tfrac 12\{v,w\}\bigr).
$$

\begin{proposition}
\label{th:group}
 For any $g\in\Sp(\V_{2\infty})$ there exists a unique up to a scalar factor
 unitary operator $W(g)$ in  $\ell^2(V^d_{\infty})$
 such that
 \begin{equation}
 a(vg)=W(g)^{-1} a(v) W(g)
 \label{eq:intertwiner}
  \end{equation}
 and 
 $$
 W(g_1) W(g_2)=c(g_1,g_2) W(g_1g_2),\qquad \text{where $c(g_1,g_2)\in\C$, $|c(g_1,g_2)|=1$.}
 $$
\end{proposition}

This statement is inside the original construction of A.~Weil (the group $V_\infty^d$
is locally compact), in Section 2 we present a proof
(this proof is used in Section 3). 



{\sc Remark.} The author does not know is 
the Weil representation of $\Sp(\V_{2\infty})$ projective or linear,
for a finite symplectic group $\Sp(2n,\F)$
it is linear, see, e.g., \cite{GH}.
\hfill $\boxtimes$

\sm

{\bf \punct Linear relations.}
Let $X$, $Y$ be linear spaces over $\F$.
A {\it linear relation} $T: X\tto Y$ is a  linear subspace in $X\oplus Y$. 
Let $T:X\tto Y$, $S:Y\tto Z$ be linear relations. Their product $ST$ is defined as the set of all 
$(x,z)\in X\oplus Z$ such that there exists $y\in Y$ satisfying $(x,y)\in T$, $(y,z)\in S$. 

\sm 

{\sc Remark.}  Let $A:X\to Y$ be a linear operator. Then its graph
$\graph(A)$ is a linear relation
$X\tto Y$. If $A:X\to Y$, $Y\to Z$ are linear operators, then
	product of their graphs 
	$$\graph(A):X\tto Y,\qquad\graph(B):Y\tto Z$$
	is the linear relation $\graph(BA):X\tto Z$.
	\hfill $\boxtimes$
	
	\sm

For a linear relation $T:X\tto Y$ we define:

\sm

$\bullet$ the {\it kernel} $\ker T$ is the intersection $T\cap (X\oplus 0)$;

\sm

$\bullet$ the {\it domain} $\dom T$ is the image of the projection of $T$ to $X\oplus 0$;

\sm

$\bullet$ the {\it image} $\im T$ is  the image of the projection of $T$ to $0\oplus Y$;

\sm

$\bullet$  the {\it indefiniteness} $\indef T$ is the intersection $T\cap (0\oplus Y)$.

\sm

We  define the {\it pseudo-inverse relation} $T^\square:Y\tto X$ as the image of $T$ under
the  transposition of summands in $X\oplus Y$.

\sm

For a linear relation $T:X\tto Y$ and a subspace $U\subset X$ we define 
a subspace $TU\subset Y$ consisting of all $y\in Y$ such that there exists 
$u\in U$ satisfying $(u,y)\in T$.

\sm

{\sc Remark.} For a linear operator $A:X\to Y$ we have
\begin{align*}
\ker\graph(A)=\ker A, \qquad\dom\graph(A)=X, \\
\im\graph A=\im(A),\qquad
\indef\graph(A)=0.
\end{align*}
 If $A$ is invertible, then $\graph(A)^\square=\graph(A^{-1})$.
For a subspace $U\subset X$ we have $\graph(A)U=AU$.
\hfill $\boxtimes$

\sm

{\bf \punct Perfect Lagrangian linear relations.}
We say that a  subspace $Z\subset \V_{2\infty}$ is {\it codiscrete}
if it is closed and the quotient $\V_{2\infty}/Z$ is discrete.

Let $\mu$, $\nu=0$, 1, \dots, $\infty$. Define a skew-symmetric bilinear form
in $\V_{2\mu}\oplus \V_{2\nu}$ by 
$$
\bigl\{(v,w),(v',w')\bigr\}_{\V_{2\mu}\oplus \V_{2\nu}}=
\bigl\{v,v'\bigr\}_{\V_{2\mu}}- \bigl\{w,w'\bigr\}_{\V_{2\nu}}
.
$$

We say that a linear relation $T: \V_{2\mu}\tto \V_{2\nu}$
is {\it perfect Lagrangian} if it 
satisfies the following conditions:

\sm

1) subspace $T$ is maximal isotropic in $\V_{2\mu}\oplus\V_{2\nu}$ (in particular, $T$ is closed);

\sm

2) $\ker T$ and $\indef T$ are compact;

\sm

3)  $\im T$, $\dom T$ are  codiscrete.

\sm

{\sc Remark.} Let $A\in \Sp(\V_{2\mu})$. Then $\graph(A)$
is a perfect Lagrangian linear relation $\V_{2\mu}\tto \V_{2\mu}$.
\hfill $\boxtimes$

\begin{lemma}
\label{l:1}
 For any perfect Lagrangian linear relation $T$,
 \begin{align*}
 (\dom T)^\lozenge=\ker T,\qquad (\ker T)^\lozenge=\dom T,
 \\
  (\im T)^\lozenge=\indef T,\qquad (\indef T)^\lozenge=\im T.
 \end{align*}
\end{lemma}

\begin{lemma}
	\label{l:2}
	Let $T:\V_{2\mu}\to \V_{2\nu}$ be a perfect Lagrangian linear relation
	and $\ker T=0$, $\indef T=0$. Then $\mu=\nu$ and $T$ is a graph
	of an element of $\Sp(\V_{2\mu})$.
\end{lemma}

\begin{theorem}
\label{th:relations}
If linear relations $T: \V_{2\mu}\tto \V_{2\nu}$, $S: \V_{2\nu}\tto \V_{2\kappa}$
are perfect Lagrangian, then $ST: \V_{2\mu}\tto \V_{2\kappa}$
is perfect Lagrangian.
\end{theorem}

We define a category $\cL$, whose objects are the spaces $\V_{2\mu}$ and morphisms
 are perfect Lagrangian relations. It is more convenient to think that objects 
 of the category are linear spaces, which are equivalent to the spaces $\V_{2\mu}$
 as linear spaces equipped with topology and skew-symmetric bilinear form.
 
 \sm

{\bf\punct Extension of the Weil representation.}

\begin{theorem}
\label{th:functor}
{\rm a)} Let $T:\V_{2\mu}\tto\V_{2\nu}$ be a perfect Lagrangian
linear  relation. Then were exists a unique up to a scalar
factor non-zero bounded operator $W(T):\ell^2(V_{\mu}^d)\to \ell^2(V_{\nu}^d)$ such that
\begin{equation}
a(w) W(T)= W(T) a(v),\qquad \text{for any $(v,w)\in T$.} 
\label{eq:def-functor}
\end{equation}

{\rm b)} For any perfect Lagrangian linear relations   $T:\V_{2\mu}\tto \V_{2\nu}$, 
$T:\V_{2\nu}\tto \V_{2\kappa}$,
$$
W(S)\,W(T)=c(T,S)\, W(ST), \qquad \text{where $c(T,S)\in \C^\times$.}
$$

{\rm c)} For any perfect Lagrangian $T:\V_{2\mu}\tto \V_{2\nu}$,
$$
W(T)^*=W(T^\square).
$$
\end{theorem}

{\bf \punct Gaussian operators.} 
Let $H\subset V^d_\mu\oplus V^d_\nu$ be a subspace such that projection
operators $H\to V^d_\nu\oplus 0$ and $H\to 0\oplus V^d_\mu$ have finite
dimensional kernels. Let $Q$ be a quadratic form on $H$.
We define a Gaussian operator $G(H;Q):\ell^2(V_\mu^d)\to \ell^2(V_\nu^d)$
by the formula
\begin{equation}
\label{eq:gauss}
G(H;Q) f(x)= \sum_{y:\,(y,x)\in H} \Exp\bigl( Q(x,y)\bigr) \,f(y)
.
\end{equation}

\begin{theorem}
 \label{th:gauss}
 Any operator $W(T):\ell^2(V_\mu^d)\to \ell^2(V_\nu^d)$ is a Gaussian operator up to a scalar factor.
 Any Gaussian operator has such a form.
\end{theorem}

{\bf\punct Bibliographic remarks.%
	\label{ss:remarks}}
1) In 1951, K.~O.~Friedrichs \cite{Fri} noticed that the Stone--von Neumann theorem
implies the existence of a representation of a real symplectic group $\Sp(2n,\R)$.
The finite-dimensional groups $\Sp(2n,\R)$ were not interesting to him and  he
formulated a conjecture  about an infinite-dimensional symplectic  group.
 The
representation of $\Sp(2n,\R)$ was explicitly described by I.~Segal in 1959,  \cite{Seg}.
A.~Weil extended this construction to groups over $p$-adic and finite fields.

\sm

2) For infinite-dimensional symplectic group the 'Weil representation' (which was conjectured
by Friedrichs) was constructed
by D.~Shale \cite{Sha} and F.~A.~Berezin \cite{Ber-dan}, \cite{Ber-book}.
The $p$-adic case was examined by M.~L.~Nazarov \cite{NNO}, \cite{Naz}
and E.~I.~Zelenov \cite{Zel}.

\sm

3) The third part of the topic was a categorization of the Weil representation,
 it had a long 
pre-history (see \cite{How}), for real and $p$-adic infinite-dimensional
groups it was  obtained by M.L.Nazarov, G.I.Olshanski, and the author
in \cite{NNO}, 1989 (see more detailed expositions of different topics
in \cite{Ner-faa}, \cite{Naz}, \cite{Olsh}, \cite{Ner-cat}).

For $\Sp(2n,\F)$ the corresponding representation of the semigroup of
linear relations was done in the preprint of R.Howe (1976),
which was not published as a paper. A detailed exposition of
this representation is contained in%
\footnote{I refered to \cite{NNO} being sure that case of finite
	fields was considered in that paper.} \cite{Ner-gauss}, Section 9.4.
 
 \sm

4) Now there exists  well-developed representation theory of infinite-dimensional
classical groups and infinite symmetric groups.
Numerous efforts  to develop a parallel representation theory of infinite-dimensional
classical groups over finite fields met some obstacles, which apparently,
were firstly observed by Olshanski in \cite{Olsh}, 1991.
For  a collection of references on different
approaches  to this topic, see \cite{Ner-finite}).
In \cite{Ner-finite} there was proposed a version of $\GL(\infty)$
over finite fields as a group of continuous operators in $\V_{2\infty}$,
it seems that this group is an interesting and hand object from
the point of view of representation theory. In particular,
there is a natural $\GL(\V_{2\infty})$-invariant measure on 
an infinite-dimensional Grassmannian (and also on flag spaces)  over $\F$ and explicit
decomposition of $L^2$ on this Grassmannian. Clearly, this topic
is related to the group $\Sp(\V_{2\infty})$ and  this  requires an 
examination of the Weil representation.

\sm

{\bf \punct Acknowledges.} This note is a kind of addition to
the work by
Nazarov, Neretin, and Olshanski \cite{NNO},
\cite{Ner-faa}, \cite{Naz}, \cite{Olsh}. I am  grateful to M.L.Nazarov and G.I.Olshanski.
I also thank R.Howe who recently sent me preprint \cite{Howe}.

\sm

{\bf\punct Structure of the paper.}
In Section 2 we prove Proposition \ref{th:group}.
The proof follows I.Segal arguments \cite{Seg},
see also \cite{Ner-gauss}, Section 1.2. The only additional
difficulty is Proposition \ref{pr:generators} about generators of the group
$\Sp(\V_{2\infty})$.


Statements on linear relations are established in Section 3.
A proof of the theorem about product of perfect Lagrangian relations
is based on all the remaining statements, this makes a structure
of the section slightly sophisticated.
 Lemmas \ref{l:1} and \ref{l:2} are proved in
Subs. \ref{ss:submodule}--\ref{ss:l:2}, Theorem \ref{th:functor}.a in
Subs. \ref{ss:th-functor-a}. Properties of Gaussian operators
(Theorem \ref{th:gauss}) are established in 
Subs. \ref{ss:gauss-bounded}--\ref{ss:gauss-end}.
Theorem \ref{th:relations} is proved
in Subs. \ref{ss:pr-proof}.
The remaining part of Theorem \ref{th:functor} is proved in
Subs. \ref{ss:end}.

\section{The Weil representation of the group $\Sp(\V_{2\infty})$.}

\COUNTERS

Here we prove Proposition \ref{th:group}.

\sm

{\bf\punct Preliminary remarks on linear transformations in $\V_{2\infty}$.\label{ss:preliminary}}
Denote by $\GL(\V_{2\infty})$ the group
of all continuous invertible linear transformations $v\mapsto vg$ of $\V_{2\infty}$.
Let 
\begin{equation}
g=\begin{pmatrix}
        a&b\\c&d
       \end{pmatrix}\in \GL(\V_{2\infty}).
       \label{eq:GL}
\end{equation}
Then the block $c$ contains only finite number of nonzero elements;
each row of $a$ contains only a finite number of nonzero elements, and each comumn of $d$
contains a finite number of nonzero elements. 

We need some statements from \cite{Ner-faa}, Sect. 2.

\sm

{\bf A}. {\it Let $g: \V_{2\infty}\to \V_{2\infty}$ be a continuous surjective
	linear transformation
	in $\V_{2\infty}$. Then the inverse transformation is continuous.}

\sm

{\bf B}. We say that a continuous linear transformation $P$ in $V^d_\infty$ or $V^c_\infty$
is {\it Fredholm} if $\dim\ker P<\infty$, $\im P$ is closed (this condition is nontrivial only for $V^d_\infty$), and $\codim \im P<\infty$.
 The {\it index}  of a Fredholm operator $P$
is 
$$\ind P:=\dim \ker P-\codim \im P.$$
The following statements hold:

\sm

1. {\it An operator $P$ in $V^c_\infty$ is Fredholm iff the dual operator 
 	$P^t$ in  $V^d_\infty$ is Fredholm}.

\sm

2. {\it Let $P$ be Fredholm and $H$ have a finite rank. Then $P+H$
	is Fredholm and }
$$
\ind (P+H)=\ind P.
$$

3. {\it If $P$, $Q$ are Fredholm linear transformations of
$V^c_\infty$ or $V^d_\infty$, then} 
$$
\ind PQ=\ind P+\ind Q.
$$

\sm

{\bf C}. {\it For any matrix $g=\begin{pmatrix}
	a&b\\c&d
	\end{pmatrix}\in\GL(\V_{2\infty})$ {\rm (see (\ref{eq:GL}))}
		the blocks $a$, $d$ are Fredholm and $\ind d=-\ind a$.}
	
	\sm

{\bf\punct Fredholm indices of blocks of $g\in\Sp(\V_{2\mu})$.}
Elements $g\in\Sp(\V_{2\mu})\subset \GL(\V_{2\mu})$ satisfy obvious additional
conditions, below we  refer to
\begin{eqnarray}
 c^t a=a c^t; \label{eq:cond-1}
 \\
 ad^t-bc^t=1. \label{eq:cond-2}
\end{eqnarray}

\begin{lemma}
	\label{l:Fredholm}
	Let $g=\begin{pmatrix}
	a&b\\c&d
	\end{pmatrix}\in \Sp(\V_{2\mu})$. Then $\ind a=\ind d=0$.
\end{lemma}

{\sc Proof.} For $g\in \Sp(\V_{2\mu})$, we have 
$$
\ind (ad^t)=\ind (a)+\ind(d^t)=\ind (a)-\ind(d) =2\ind(a).
$$
On the other hand, $bc^t$ has finite rank. By (\ref{eq:cond-2}),
$$
\ind (ad^t)=\ind (1+bc^t)=\ind (1)=0.
$$
Thus, $\ind (a)=0$.
\hfill $\square$

\sm

{\bf \punct Generators of the group $\Sp(\V_{2\infty})$.\label{ss:preliminary-s}}
We define the following subgroups in $\Sp(\V_{2\infty})$:

\sm

---  the subgroup $H$
consists of matrices of the form
$\begin{pmatrix}
      a&0\\0&a^{t-1}
\end{pmatrix}$;

\sm

--- the subgroup
$N_+$ consists of matrices of the form $\begin{pmatrix}
                             1&b\\0&1
                            \end{pmatrix}$, where $b=b^t$;
                            
\sm                            
                             
--- the subgroup $N_-$ consists of matrices of the form 
$\begin{pmatrix}
                             1&0\\c&1
                            \end{pmatrix}$, where $c=c^t$.
                            
                            \sm
                            
  Denote by $J_k$ the following block matrix of the size
$((k-1)+1+\infty)+((k-1)+1+\infty)$:
$$J_k=
\begin{pmatrix}
1&0&0&0&0&0 \\
0&0&0&0&1&0\\
0&0&1&0&0&0\\
0&0&0&1&0&0\\
0&-1&0&0&0&0\\
0&0&0&0&0&1\\
\end{pmatrix}  
.
$$

  \begin{proposition}
  	\label{pr:generators}
    	The group $\Sp(\V_{2\infty})$ is generated by the subgroup $H$, $N_+$, and elements $J_k$.
                                  \end{proposition}

{\bf\punct Generators of $\Sp(\V_{2\infty})$. Proof of Proposition \ref{pr:generators}.}      Let $s=\begin{pmatrix}1&0\\c&1\end{pmatrix}\in N_-$.
Since $Q$ contains only a finite number of nonzero elements, we have
$$
J_1\dots J_m s J_m^{-1}\dots J_1^{-1}\in N_+
$$
for sufficiently large $m$.
Thus,   $N_-$ is contained  in the subgroup generated by $H$, $N_+$, and $J_k$. 

Let
$$g= \begin{pmatrix}
        a&b\\c&d
       \end{pmatrix}\in \Sp(\V_{2\infty}).$$
       We will multiply it by elements of groups $H$, $N_+$, $N_-$ and elements $J_k$ in an appropiate  way.
       As a result we will come to the unit matrix.
       
       \sm
       
  Recall that $a$ is Fredholm of index 0. Therefore  there are invertible operators $K$, $L$
  in $V_\infty^d$ such that 
  $
  K a L
  $
  is a block $(l+\infty)\times (l+\infty)$-matrix of the form $  \begin{pmatrix}
   0&0\\0&1
  \end{pmatrix}$. This follows from \cite{Ner-finite}, Lemma 2.7.
  We pass to a new matrix $g'$ given by
  $$
g':=  \begin{pmatrix}
   K&0\\0&K^{t-1}
  \end{pmatrix}
\begin{pmatrix}
        a&b\\c&d
       \end{pmatrix}
    \begin{pmatrix}
   L&0\\0&L^{t-1}
  \end{pmatrix}=\begin{pmatrix}
  0&0&b'_{11}&b'_{12}\\
  0&1&b'_{21}&b'_{22}\\
  c'_{11}&c'_{12}&d'_{11}&d'_{11}\\
  c'_{21}&c'_{22}&d'_{21}&d'_{22}
  \end{pmatrix}  
  $$
      (the size of the right hand side is $l+\infty+l+\infty$). The condition
      (\ref{eq:cond-1}) implies $c'_{12}=0$, $c'_{21}=0$, $(c'_{22})^t=c'_{22}$.
      In particular the following matrix
      $$
s:=\begin{pmatrix}
    1&0&0&0\\
    0&1&0&0\\
    0&0&1&0\\
    0&0&-c'_{22}  &1 
   \end{pmatrix} 
        $$
        is contained in $N_-$.
The element $s g'$ has the form
$$
sg'=\begin{pmatrix}
  0&0&b''_{11}&b''_{12}\\
  0&1&b''_{21}&b''_{22}\\
  c''_{11}&0&d''_{11}&d''_{11}\\
  0&0&d''_{21}&d''_{22}
  \end{pmatrix}=:g''.
$$
The matrix $ c''_{11}$ is nondegenerate
(otherwise, the whole matrix $sg'$ is degenerate). We take an element
$$
J_1\dots J_l g''=:g'''
$$
and get a matrix of the form
$$
g'''=
\begin{pmatrix}
 -c''_{11}&0&b'''_{11}&b'''_{12}\\
 0&1&b'''_{21}&b'''_{22}\\
 0&0&d'''_{11}&d'''_{12}\\
  0&0&d'''_{21}&d'''_{12}
\end{pmatrix}.
$$
Keeping in the mind (\ref{eq:cond-2}) we observe
$$
\begin{pmatrix}
 d'''_{11}&d'''_{12}\\
 d'''_{21}&d'''_{12}\\
\end{pmatrix}
=
\begin{pmatrix}
 -c''_{11}&0\\
  0&1
  \end{pmatrix}^{-1}.
$$
Element
$$
\begin{pmatrix}
-c_{11}'''&0&0&0\\
0&1&0&0\\
0&0& -(c_{11}''')^{-1}&0\\
0&0&0&1
\end{pmatrix}^{-1} g'''
$$
is contained in $N_+$.
\hfill $\square$

\sm

{\bf \punct Construction of the Weil representation of the symplectic groups $\Sp(\V_{2\infty})$.%
\label{ss:construction-group}}
A remaining part of a   proof of Proposition \ref{th:group} is based on standard  arguments (see \cite{Ner-gauss}, Sect. 1.2), which were
proposed by I.~Segal).

\sm

{\sc Step 1.} {\it The representation $a(v)$ of the Heisenberg group in $\ell^2(V_\infty^d)$} is
irreducible. Indeed, let us show that there are no nontrivial
intertwining  operators. 
Any bounded operator $Q$ in $\ell^2(V_\infty^d)$
commuting with all  operators
$$
a(0,v^c) f(x)=\Exp\Bigl(\sum v^c_j x_j\Bigr)\,f(x)
$$
is an operator of multiplication by a function.
Since $Q$ commutes also with shifts
$$
a(v^d,0) f(x)=f(x+v^d),
$$
we get that $Q$ is a multiplication by a constant. 

\sm

{\sc Step 2.}
{\it An operator $W(g)$ is defined up to a scalar factor {\rm(}if it exists}).
Indeed, the map $v\mapsto vg$ is an automorphism of the Heisenberg group.
Therefore, the formula $v\mapsto a(vg)$ determines
a unitary representation of the Heisenberg group. The operator
$W(g)$ intertwines unitary representations $a(v)$ and $a(vg)$. By the Schur
lemma, $W(g)$ is unique up to a scalar factor.

\sm

{\sc Step 3.} If $W(g_1)$, $W(g_2)$ exist, then
$$
W(g_1)W(g_2)=\lambda\cdot W(g_1 g_2).
$$
Indeed,
\begin{multline*}
\bigl(W(g_1) W(g_2)\bigr)^{-1}a(v)W(g_1) W(g_2)=
W(g_2)^{-1}\bigl( W(g_1)^{-1}a(v)W(g_1)\bigr) W(g_2)=\\=
W(g_2)^{-1} a(vg_1) W(g_2)=a(vg_1g_2).
\end{multline*}

{\sc Step 4.} It remains to write operators corresponding
to generators of the group $\Sp(\V_{2\mu})$.

\sm

--- Operators corresponding to elements of the subgroup $H$
are 
$$
W\begin{pmatrix}
    a&0\\0&a^{t-1}
   \end{pmatrix} f(x)=f(xa);
$$

--- Operators corresponding to elements of $N_+$ are
\begin{equation}
W\begin{pmatrix}
  1&b\\0&1
 \end{pmatrix} f(x)=\Exp\bigl(\tfrac 12 \sum b_{kl} x_k x_l\bigr)\,f(x).
 \label{eq:multiplication}
\end{equation}

--- An operator corresponding to $J_k$ is the Fourier transform 
with respect to the coordinate $x_k$.

\section{The Weil representation of the category of Lagrangian relations}

\COUNTERS

{\bf\punct On a canonical form of a compact isotropic submodule.
	\label{ss:submodule}}

\begin{lemma}
	\label{l:canonical}
 For any compact isotropic subspace $M\subset\V_{2\infty}$ there is
 an element $g\in\Sp(\V_{2\infty})$ such that $Mg\subset V_\infty^c$.
 Moreover, we can choose $g$ is such a way that $Mg\subset V_\infty^c$ be a subspace given by a system
 of equations of the type $y_\alpha=0$, where $\alpha$ ranges in a subset $A\subset\N$.
\end{lemma}

{\sc Proof.}
Consider the projection map $\pi:M\to V^c_\infty$. Its fibers are compact
and discrete,  therefore they are finite.
In particular, $\pi^{-1}(0)$ is a finite-dimensional subspace in
$V_\infty^d$.
We can choose an element $h$ of the subgroup $H\subset\Sp(\V_{2\infty})$
  such that 
$(\pi^{-1}(0)) h\subset V^d_\infty$ is the subspace consisting of vectors $(x_1,\dots,x_k,0,\dots)$.
Since $M h$ is isotropic, $\pi(Mh \bigr)$
is contained in the subspace $y_1=\dots=y_k=0$. Therefore
$Mh J_1\dots J_k\subset V_\infty^c$. 

Next, we wish to describe closed subspaces in $V^c_\infty$ modulo linear transformations of $V^c_\infty$
(they are induced by elements of the subgroup $H$).
 The Pontryagin duality determines a one-to-one correspondence between
 sets of closed subgroups in Abelian groups $V^c_\infty$ and $V^d_\infty$
 (see \cite{Mor}, Theorem 27); the both groups are equipped with an Abelian group of automorphisms $x\mapsto\lambda x$, where
 $\lambda\in\F^\times$. It is easy to see that the correspondence send
 invariant subgroups (subspaces) to invariant subgroups. 
Therefore the question is reduced
 to a description of subspaces in $V^d_\infty$
 modulo linear transformations,
the latter problem is trivial.
\hfill $\square$

\begin{corollary}
 Let $L$ be a compact isotropic submodule of $\V_{2\infty}$.
 Then the space $L^\lozenge/L$ is isomorphic as a symplectic space to
 $\V_{2n}$ or $\V_{2\infty}$.
\end{corollary}



\begin{lemma}
 For any codiscrete coisotropic submodule $L\subset \V_{2\infty}$
 there exists $g\in\Sp(\V_{2\infty})$ such that $Lg\supset V_\infty^c$.
\end{lemma}

{\sc Proof.} We reduce the compact isotropic module $L^\lozenge$
to the canonical form (i.e., $L^\lozenge$ is a coordinate subspace in $V_\infty^c$).
\hfill $\square$

\sm

{\bf\punct Proof of Lemma \ref{l:1}.\label{ss:l:1}}
 Consider a perfect Lagrangian
relation $T$. Obviously, for $v\in\ker T$, $Z\in \dom T$, we have
$\{v,z\}=0$, i.e., $(\dom T)^\lozenge \supset \ker T$.
 Next, $(\dom T)^\lozenge = \ker T$; otherwise,
we take an isotropic linear relation $T+(\dom T)^\lozenge\supsetneq T$.
\hfill $\square$

\sm

{\bf\punct Proof of Lemma \ref{l:2}.\label{ss:l:2}}
Only the case $\mu=\nu=\infty$ requires a proof.
Let $T:\V_{2\infty}\tto \V_{2\infty}$ be a perfect Lagrangian linear relation,
$\ker T=0$, $\indef T=0$. According Lemma \ref{l:1}, we
have
$\dom T=\V_{2\infty}$, $\im T=\V_{2\infty}$.

\begin{lemma}
	\label{l:graph-topology}
	The subspace
$T$  is
isomorphic to $\V_{2\infty}$ as a linear space with a topology.
\end{lemma} 

{\sc Proof of Lemma \ref{l:2}.} We apply statement {\bf A} of
Subsect. \ref{ss:preliminary}. The projections
$$\pi_T:T\to \V_{2\infty}\oplus 0, \qquad \pi_T':T\to 0\oplus\V_{2\infty}$$
 are continuous, therefore the inverse map
$\pi_T^{-1}$ is continuous. Hence the map $\pi_T'\circ\pi_T^{-1}$
is continuous, this is the the linear transformation whose graph is $T$.
\hfill $\square$

\sm

{\sc Proof of Lemma \ref{l:graph-topology}.}
Represent $\V_{2\infty}$ as a union of a chain of compact subspaces
$$
W_0=V_\infty^c\subset W_1\subset W_2\subset\dots
$$
Consider the projection map
 $\pi:\V_{2\infty}\oplus \V_{2\infty}\to \V_{2\infty}\oplus 0$,
 denote by $\pi_T$ its restriction  to $T$. Then
 $$
 \pi_T^{-1} V^c_\infty=\cup_{j=1}^\infty \bigl[T\cap (V^c_\infty\oplus W_j)\bigr]
 .
 $$
 Therefore
 $$
   V^c_\infty=\cup_{j=1}^\infty \pi_T \bigl[T\cap (V^c_\infty\oplus W_j)\bigr]
 .$$
 In the right hand side we have a union of an increasing sequence of compact sets,
 the left hand side is compact. Therefore for some $k$,
 $$
 V^c_\infty= \pi_T \bigl[T\cap (V^c_\infty\oplus W_k)\bigr].
 $$
 The set $T\cap (V^c_\infty\oplus W_k)$ is compact and  $\pi_T$
 is continuous, therefore the inverse map $\pi_T^{-1}:V^c_\infty
 \to T\cap (V^c_\infty\oplus W_k)$ is continuous.
 
 
Next, denote by $e_l$  the standard basis in $ V^d_\infty$.
Then 
$$
T\simeq \bigl[T\cap (V^c_\infty\oplus W_k)\bigr]\oplus 
\bigoplus_{l=1}^\infty \F\cdot \pi_T^{-1} e_l.
$$
Thus $T$ is isomorphic $\V_{2\infty}\oplus0$. 
\hfill $\square$


\sm
 
{\bf\punct Proof of Theorem \ref{th:functor}.a. Existence
	of operators $W(T)$.%
	\label{ss:th-functor-a}}

\begin{lemma}
 Let the claim of Theorem {\rm \ref{th:functor}.a} hold for a linear relation
$T:\V_{2\mu}\tto \V_{2\nu}$. Then the same statement holds for any linear relation
$gTh$, where $g\in\Sp(\V_{2\nu})$, $h\in \Sp(\V_{2\mu})$. 
\end{lemma}

{\sc Proof.} Obvious.
\hfill $\square$

\sm

We start a proof of the theorem.
Keeping in the mind Lemma \ref{l:canonical} we can assume
$$
\ker T\subset V_\mu^c,\qquad \indef T\subset V_\nu^c.
$$
Let $w\in\indef T$. Then
$$a(w)\,W(T)= W(T)a(0)=
W(T).$$
The operator $a(w)$ is an operator
of multiplication by a function, 
$$a(w)f(x)= \Exp \Bigl(\sum x_j w_j^c\Bigr) f(x)$$
Therefore any function $\psi\in \im W(T)$
vanishes on the set $\Exp \bigl(\sum x_j w_j^c\bigr)\ne 1$.
 In other words, all elements of $\psi\in\im W(T)$
 are supported by the subspace
  $(\indef T)^\circ \subset V_\nu^d$.

Let $v\in \ker T$. Then
$$W(T)=W(T)a(v).$$
The operator  $a(v)$ is a multiplication by function taking finite number of values
 $\lambda_l=\exp(2\pi li/p)$. Therefore, $\lambda_l$ are the eigenvalues of $a(v)$,
 and $T(W)=0$ on all subspaces $\ker(a(v)-\lambda_l)$ for $\lambda_l\ne 1$
Therefore, $W(T)$ is zero for any function supported by $V_\mu^d\setminus (\ker T)^\circ$.

Consider decompositions
\begin{align*}
\ell_2(V^d_\mu)&=\ell_2\bigl((\ker T)^\circ\bigr)\oplus
 \ell_2\bigl(V^d_\mu\setminus(\ker T)^\circ\bigr);
 \\
\ell_2(V^d_\nu)&=\ell_2\bigl((\indef T)^\circ\bigr)\oplus
\ell_2\bigl(V^d_\nu\setminus(\indef T)^\circ\bigr). 
\end{align*}
The operator 
$W(T)$ has the following block form with respect to this decomposition
$$
W(T)=\begin{pmatrix}
\wt W(T)&0\\0&0
\end{pmatrix}
,$$
with a non-zero block 
$$
\wt W(T) :\,\ell^2\bigl((\ker T)^\circ\bigr)\,\to\, \ell^2\bigl((\indef T)^\circ\bigr).
$$

\sm

The linear relation $T$ determines a linear relation $T':\dom T\tto\im T$.
We take the projection
$$
\dom T\oplus\im T\,\to \,(\dom T/\ker T)\oplus(\im T/\indef T)
$$
and the image $\wt T$ of $T'$ under this projection.
Thus we get a linear relation $\wt T:\dom T/\ker T\tto \im T/\indef T$.
The spaces $\dom T/\ker T$, $\im T/\indef T$ have form
$\V_{2\kappa}$. By construction, $\ker \wt T=0$, $\indef \wt T=0$.
Therefore, $\wt T$ is a graph of a symplectic operator in $\V_{2\kappa}$.

It is easy to that 
$\wt W(T)$
has the same commutation relations with the Heisenberg group as 
$W(\wt T)$. It remains to refer to Proposition \ref{th:group}.
\hfill $\square$



\sm

{\bf \punct One corollary from the previous Subsection.}
For a linear embedding $B:V^d_\mu\to V^d_\nu$ we define two 
operators
$$
\sigma_B: \ell^2(V^d_\nu)\to \ell^2(V^d_\mu),\qquad
\sigma_B^*: \ell^2(V^d_\mu)\to \ell^2(V^d_\nu)
$$
in the following way 
\begin{align*}
 \sigma_B \phi(x)&= \phi(xB);
 \\
 \sigma^*_B \psi(y)&=\begin{cases}
              \psi\bigl(B^{-1}(y)\bigr),\qquad &\text{if $y\in \im B$};\\
              0,\qquad &\text{if $y\notin\im B$}.
                    \end{cases}
\end{align*}

In fact, in the previous section we proved the following lemma:

\begin{lemma}
\label{l:WT-decompose}
 Any operator $W(T)$ can be decomposed as product of the form
 $$
 W(T)=W(g_1) \lambda^*_B W(g_2) \lambda_C W(g_3),
 $$
 where $g_1$, $g_2$, $g_3$ are elements of symplectic groups and $B$, $C$
 are appropriate embeddings.
\end{lemma}

{\bf \punct Gaussian operators are bounded.%
	\label{ss:gauss-bounded}}
Consider a Gaussain operator $G(H,Q)$ given by (\ref{eq:gauss}).
 Consider projections 
$\pi_1:H\to V^d_\mu\oplus 0$, $\pi_2:H\to 0\oplus V_\nu^d.$ 
We represent $H$ as a finite disjoint  union of affine subspaces
$Z_j$ such that  $\pi_1$, $\pi_2$ are injective on $Z_j$. 
Consider operators 
$$
G_j f(x)= \sum_{y:\,(y,x)\in Z_j} \Exp\bigl(  Q(x,y)\bigr) \,f(y)
.
$$
Actually, for each $x$ the sum consists of 0 or 1 element. Clearly, $\|G_j\|=1$,
and $G(H,Q)=\sum G_j$.

\sm

{\bf \punct Products of Gaussian operators.}

\begin{lemma}
\label{l:gauss-sums}
 Let $X$ be a finite-dimensional space, $Y$ be a sum of a finite or countable 
 number of copies of $\F$. Let $Q$ be a quadratic form
 on $X\times Y$. Consider the sum
 \begin{equation}
 F(y)=\sum_{x\in X} \Exp\bigl(Q(x,y)\bigr).
 \label{eq:gauss-sum}
 \end{equation}
 Then there is a subspace $Z\subset Y$ of  codimension $\le\dim X$, a nonzero constant $c$,
  and a quadratic form $R$ on $Z$ such that
 $$
 F(y)=\begin{cases}
       c\cdot \Exp\bigl( R(y)\bigr), \qquad &\text{\rm if $y\in Z$};
       \\
       0,\qquad &\text{\rm if $ y\notin Z$}.
      \end{cases}
 $$
\end{lemma}

{\sc Proof.} This follows from the following  observation (see, e.g., \cite{Ner-gauss}, Sect. 9.2).
Let $X$, $U$ be  finite-dimensional spaces
over $\F$ of the same dimension.  Consider a quadratic form
$P(x)$ on $X$.
Let 
$$
f(u)=
\sum_{x\in X} \Exp\bigl(P(x)+\sum u_j y_j\bigr).
$$
Then there is a subspace $K\subset U$ and a quadratic form $S$ on $U$ such that
\begin{equation}
f(u)=
\begin{cases}
       c\cdot \Exp\bigl( S(u)\bigr), \qquad &\text{if $u\in K$};
       \\
       0,\qquad &\text {if $u\notin K$}
      \end{cases}.  
      \label{eq:gauss-sums}
\end{equation}

We represent $Q(x,y)$ as
$$
Q(x,y)=Q(x,0)+\sum_j x_j l_j(y)+Q(0,y),
$$
where $l_j$ are linear forms on $Y$, and apply formula (\ref{eq:gauss-sums}).
\hfill $\square$

\begin{lemma}
 A product of Gaussian operators is a Gaussian operator up to a nonzero scalar factor.
\end{lemma}

{\sc Proof.}
Consider Gaussian operators
$$
G(H,Q):\ell^2(V_\mu^d)\to \ell^2(V_\nu^d),\qquad G(K,R):\ell^2(V_\nu^d)\to \ell^2(V_\kappa^d).
$$
We consider the set $Z$ of all triples 
$(u,v,w)\in V_\mu^d\oplus V_\nu^d\oplus V_\kappa^d$, such that
$(u,v)\in H$, $(v,w)\in K$.
The kernel of the product is given by
$$
N(u,w)=
\sum_{v:\,(u,v,w)\in Z}
\Exp\bigl(Q(u,v)\bigr) \Exp\bigl(R(v,w)\bigr).
$$
We get a sum of the form (\ref{eq:gauss-sum}).
 More precisely, consider the natural projection 
$Z\to V_\mu^d\oplus V_\kappa^d$. Denote by $X$ its kernel (it is finite-dimensional),
let $Y$ be a subspace complementary to the $X$. We apply Lemma \ref{l:gauss-sums}
and obtain a Gaussian expression for the kernel $N$.
\hfill$\square$

\begin{corollary}
 For any $g\in\Sp(\V_{2\mu})$ operators $W(g)$ are Gaussian.
\end{corollary}

{\sc Proof.} Indeed, for generators of $\Sp(\V_{2\mu})$ the operators $W(g)$ are Gaussian, see 
Subs. \ref{ss:construction-group}. Their products are Gaussian. \hfill $\square$

\begin{lemma}
 For any perfect Lagrangian linear relation $T$, the operator $W(T)$ is Gaussian.
\end{lemma}

{\sc Proof.} We refer to Lemma \ref{l:WT-decompose}.
\hfill $\square$

\sm

{\bf\punct End of proof Theorem \ref{th:gauss}.%
	\label{ss:gauss-end}}

\begin{lemma}
 Any Gaussian operator has the form $W(T)$.
\end{lemma}

{\sc Proof.} Consider a Gaussian operator $G(H;Q):\ell^2(V_{\mu}^d)\to \ell^2(V_{\nu}^d)$.
Extend the quadratic form $Q$ to $V_{\mu}^d\oplus V_{\nu}^d$ in an arbitrary way.
Represent $Q$ as
$$
Q(x,y)=Q(x,0)+Q(0,y)+\sum s_{kl} x_k y_l.
$$

Let $C(x)$ be a quadratic form on a space $V_\kappa^d$.
Denote by $G[[C]]$ the operator in $\ell^2(V_\kappa^d)$ 
given by
$$
G[[C]]\phi(x)= \Exp\bigl(C(x)\bigr)
\phi(x)
$$
Recall (see Subs. \ref{ss:construction-group}, formula (\ref{eq:multiplication})),
	 that such operators
have a form $W(g)$ for certain $g\in\Sp(\V_{2\kappa})$.
Consider a Gaussian operator
$$
\cG:=G[[-Q(x,0)]]\, G[H;Q]\, G[[-Q(0,y)]].
$$
Clearly, the statements {\it 'the operator $G(H,Q)$ has a form $W(\cdot)$'} and  {\it 'the operator $\cG$ has a form $W(\cdot)$'}
are equivalent. The operator $\cG$ has a form
$$
\cG \psi(x)=\sum_{y:\,(x,y)\in H} \Exp\Bigl(\sum s_{kl} x_k y_l\Bigr)\,\psi(y)=:
\sum_{y:\,(y,x)\in H} \Exp(xSy^t)\,\psi(y).
$$
Let us describe the set $T$ of all
 $(v,w)\in \V_{2\mu}\oplus \V_{2\nu}$ such that
$$
a(w)\cG=\cG a(v).
$$
Consider $(p,q)\in H$. Then 
\begin{equation}
\bigl((p,qS); (q,-pS)\bigr)\in T.
\label{eq:pq}
\end{equation}
Next, let $(\xi,\eta)\in  V^c_\mu\oplus V^c_\nu$
be contained in $H^\circ$. Then 
\begin{equation}
\bigl( (0,\xi);(0,\eta)\bigr)\in T.
\label{eq:xi-eta}
\end{equation}
It is easy to see that elements of forms (\ref{eq:pq}) and
(\ref{eq:xi-eta}) generate
 a perfect Lagrangian relation.
\hfill $\square$

\sm

{\bf\punct Products of perfect Lagrangian relations. Proof of Theorem \ref{th:relations}.%
	\label{ss:pr-proof}}

\begin{lemma}
Fix a perfect Lagrangian relation $T$.
Let $(\beta,\gamma)$ satisfy
$$
a(\gamma)W(T)=W(T) a(\beta).
$$
Then $(\beta,\gamma)\in T$.
\end{lemma}

{\sc Proof.}
Let $p$, $q\in\V_{2\infty}$. Then
$$
a(p)a(q)a(-p)a(-q)=\Exp\bigl(\{p,q\}\bigr).
$$
Let $(\beta,\gamma)$ does not contained in  $T$. 
Choose a vector $(u,v)\in T$ that is not orthogonal to
$(\beta,\gamma)$. This means that
$\{\beta,u\}\ne \{\gamma,v\}$.
For a field $\F$ of a prime order this implies
\begin{equation}
\Exp\bigl(\{\beta,u\}\bigr)\ne \Exp\bigl(\{\gamma,v\}\bigr).
\label{eq:ne}
\end{equation}
For an arbitrary finite field we can multiply
$(u,v)$ by a constant factor, in this way we can achieve
(\ref{eq:ne}). Next, we have
$$
a(\beta)a(\gamma)a(-\beta)a(-\gamma)
W(T)
=W(T)
a(u)a(v)a(-u)a(-v),
$$
or
$$
\Exp\bigl(\{\beta,u\}\bigr)\cdot W(T)=\Exp\bigl(\{\gamma,v\}\bigr) \cdot W(T)
$$
Hence $W(T)=0$, and we came to a contradiction. 
\hfill $\square$

\begin{lemma}
	\label{l:ST}
Let $T:\V_{2\mu}\tto \V_{2\nu}$, $S:\V_{2\nu}\tto \V_{2\kappa}$ be perfect Lagrangian 
linear relations.
 Then the linear relation $ST$ satisfies the following properties:
 
 \sm
 
  {\rm a)} The linear relation $ST$ is isotropic
  and for any $(v,w)\in ST$ we have 
  \begin{equation}
a(w)\, W(S)\,W(T)=W(S)\,W(T)\, a(v).
  \end{equation}
  
  \sm
  
  {\rm b)} $\ker ST$, $\indef ST$ are compact, and $\dom ST$, $\im ST$ are codiscrete.
  
  \sm
  
   {\rm c)} $ST$ is contained in a certain perfect Lagrangian relation $R$ such that $W(R)=W(S)W(T)$.
\end{lemma}

{\sc Proof.} 
a) Let $(v,w)$, $(v',w')\in ST$. Choose $u$, $u'\in\V_{2\nu}$ such that
$(v,u)$, $(v',u')\in T$ and $(u,w)$, $(u',w')\in S$.
Since $T$, $S$ are isotropic, we have
$$
\{u,u'\}=\{v,v'\}=\{w,w'\}.
$$
Thus $ST$ is isotropic.

 Next,
$$
W(S) W(T) a(v)=W(S)a(u) W(T)=a(w) W(S) W(T)
.
$$

b) The subspace $\ker ST$ is the set of all $v$ such that there is $u$ satisfying $u\in\ker T$,
$(v,u)\in S$. Thus we take the preimage $Z$ of $\ker T$ under the projection
$T\to  V_\mu \oplus 0$, and $\ker ST$ is the projection of $Z$ to $0\oplus V_\mu$.
Fibers of the projection $T\to  V_\mu \oplus 0$ are compact, therefore $Z$
is compact, and  $\ker ST$ is compact.

 Let us verify the statement about $\im ST$.
 The subspace $\im T\cap\dom S$ is codiscrete
  in $\dom S$. Its image $H$ under the projection
  $\dom S\to \dom S/\ker S$ also is codiscrete.
  The relation $S$ determines a symplectic isomorphism
  $\wt S: \dom S/\ker S\to \im S/\indef S$.
  The subspace $\wt S H$ is codiscrete in $\im S/\indef S$.
  Therefore its lift to $\im S$ (it is $\im ST$) is codiscrete
  in $\im S$ and therefore in $\V_{2\kappa}$.

\sm

c) Operators $W(T)$ and $W(S)$ are Gaussian. Therefore, $W(S) W(T)$
is Gaussian, and hence it has a form $W(R)$.
\hfill$\square$

\begin{lemma}
	\label{l:last0}
 Let $X$, $Y$ be codiscrete coisotropic subspaces in
 $\V_{2\infty}$. Consider the symplectic space
 $X/X^\lozenge$, the image $K$ of $X\cap Y\subset X$ in 
 $X/X^\lozenge$ and the image $L$ of $X\cap Y^\lozenge\subset X$
 in  $X/X^\lozenge$. Then $L=K^\lozenge$ in  $X/X^\lozenge$.
 Moreover, $L$ is compact, and $K$ is codiscrete.
\end{lemma}

{\sc Proof.}
We have
\begin{align*}
K=\wt K/X^\lozenge, \quad\text{where} \quad
\wt K= (X\cap Y)+X^\lozenge;
\\
L=\wt L/X^\lozenge, \quad\text{where} \quad
\wt L= (X\cap Y^\lozenge)+X^\lozenge.
\end{align*}
The space $X$ is equipped with a  skew-symmetric bilinear form,
whose kernel is $X^\lozenge$. It is clear that $\wt K$ and $\wt L$
are orthogonal in $X$, therefore $K$ and $L$ are orthogonal in
$X/X^\lozenge$. Next,
$$
\wt K^\lozenge=(X^\lozenge+Y^\lozenge)\cap X.
$$
Let $h\in \wt K^\lozenge$. Let $\wt h$ be its representative, $\wt h= a+b$, where
$a\in X^\lozenge$, $b\in Y^\lozenge$. Then $b$ also is a representative of $h$,
and hence $h\in L$.
\hfill$\square$

\sm

\begin{lemma}
	\label{l:last}
	Let $T:\V_{2\mu}\tto \V_{2\nu}$ be a perfect Lagrangian linear relation.
	Let $R\subset T$ be a relation $\V_{2\mu}\tto \V_{2\nu}$. Assume that
$$(\ker R)^\lozenge=\dom R,\qquad (\dom R)^\lozenge= \ker R,$$
or
$$
(\im R)^\lozenge=\indef R,\qquad (\indef R)^\lozenge= \im R.
$$
Then $T=R$.
\end{lemma}

{\sc Proof.} Let $R\subsetneq T$.
 Considering projection to $\V_{2\mu}\oplus 0$,
we get 
$$
\dom R\subsetneq \dom T,\qquad\text{or}\qquad \indef R\subsetneq \indef T
.
$$
Therefore,
$$
\dom R\subsetneq (\ker T)^\lozenge\subset (\ker R)^\lozenge
,
$$
or
$$
\qquad \indef R\subsetneq \indef T=(\im T)^\lozenge\subset 
(\im R)^\lozenge.
$$
This contradicts to the conditions of the lemma.
\hfill $\square$, 

\sm

{\sc Proof of Theorem \ref{th:relations}.}
Let $T:\V_{2\mu}\tto \V_{2\nu}$, $S:\V_{2\nu}\tto \V_{2\kappa}$
be perfect Lagrangian linear relations. We wish to prove that 
$ST$ is perfect Lagrangian.

Without loss of generality we can assume 
\begin{equation}
\ker T=0,\qquad\indef S=0.
\label{eq:ker-T-indef-S}
\end{equation}
Otherwise, we  take the natural projection 
$$
\dom T\oplus \V_{2\nu}\,\to\, (\dom T/\ker T)\oplus \V_{2\nu}
$$
and the image $\wt T$ of $T$ under this projection.
We get a linear relation $\wt T:\dom T/\ker T\tto \V_{2\nu}$.
Clearly, the statements '{\it $ST$ is perfect Lagrangian}'
and '{\it $S\wt T$ is perfect Lagrangian}' are equivalent.
In the same way we can assume $\indef S=0$.

Under the condition (\ref{eq:ker-T-indef-S}),
$$\dom T= \V_{2\mu}, \qquad \im S=\V_{2\kappa}.$$
Next,
\begin{align*}
\im ST&= S\im T=S(\im T\cap \dom S) , \\
 \indef ST&= S \indef T= S(\indef T\cap \dom S)
 .
\end{align*}
Consider the  map $\dom S\oplus V_{2\kappa}\to (\dom S/\ker S)\oplus V_{2\kappa}$.
Let $\wh S$ be the image of $S$ under this map. It is a graph 
of an symplectic bijection $\sigma:\dom S/\ker S\to V_{2\kappa}$.
We have
\begin{align*}
\im ST&= \sigma \Bigl[\bigl((\im T\cap \dom S)+\ker S\bigr)/\ker S\Bigr],\\
\indef ST&=\sigma \Bigl[ \bigl((\indef T\cap \dom S)+\ker S\bigr)/\ker S \Bigr]
.
\end{align*}
By Lemma \ref{l:last0}, the spaces in the square brackets are orthogonal complements one to
another. Therefore, the same holds for the left hand sides.
By Lemma \ref{l:last}, $ST$ is Lagrangian.
\hfill $\square$

\sm

{\bf\punct End of proof Theorem \ref{th:functor}.%
	\label{ss:end}}
After the establishment of Theorem \ref{th:relations} Lemma \ref{l:ST}.a
becomes the statement b) of Theorem \ref{th:functor}.

To prove the statement c) of Theorem \ref{th:functor}
we write adjoint operators to  the both sides of (\ref{eq:def-functor}),
$$
W(P)^* a(-w)=a(-v) W(P)^*,
$$
or
$$
a(v)W(P)^*= W(P)^* a(w).
$$
This is the defining relation for the operator $W(P^\square)$.

	\noindent
	\tt Math.Dept., University of Vienna,
	\\
	Oskar-Morgenstern-Platz 1, 1090 Wien;
	\\
	\& Institute for Theoretical and Experimental Physics (Moscow);
	\\
	\& Mech.Math.Dept., Moscow State University;
	\\
\&	Institute for information transmission problems (Moscow);
	\\
	e-mail: neretin(at) mccme.ru
	\\
	URL:www.mat.univie.ac.at/$\sim$neretin

\end{document}